\documentstyle[leqno,12pt]{article}

\newcommand{\thenewsection}{\bf{\thesection}\rule{1em}{0em}}
\newenvironment{newsection}{\begin{flushleft}\rule{0mm}{3ex}\thenewsection
\bf}{\rule{0mm}{3ex}\noindent\end{flushleft}}

\newcommand{\mysection}[1]{\refstepcounter{section}
\begin{newsection}{#1}\end{newsection}}

\newcommand{\thenewsubsection}{\bf{\thesubsection}\rule{1em}{0em}}
\newenvironment{newsubsection}{\begin{flushleft}\rule{0mm}{3ex}\thenewsubsection
\bf}{\rule{0mm}{3ex}\noindent\end{flushleft}}

\newcommand{\proof}{{\em Proof.\ }}

\newcommand{\qed}{\hspace*{\fill}\rule{2mm}{2mm}\vspace{5mm}}

\newcommand{\norm}[1]{\mbox{$\left\|{#1}\right\|\,$}}

\newcommand{\ra}{\rightarrow}

\newcommand{\arr}[1]{\stackrel{#1}{\longrightarrow}}

\newcommand{\uar}[1]{{\phantom{}\Big\uparrow
\makebox[0pt][r]{$\scriptstyle #1\;$}}}

\newcommand{\dar}[1]{{\phantom{}\Big\downarrow
\makebox[0pt][l]{$\scriptstyle #1$}}}

\newcommand{\Bdal}[1]{{\makebox[0pt][r]{$\scriptstyle #1$}
\phantom{}\Bigg\downarrow}}

\newcommand{\larray}{\left(\begin{array}{cc}}
\newcommand{\rarray}{\end{array}\right)}

\newcommand{\beq}{\begin{equation}}
\newcommand{\eeq}{\end{equation}}
\newcommand{\barr}{\begin{array}}
\newcommand{\earr}{\end{array}}
\newcommand{\beqar}{\begin{eqnarray}}
\newcommand{\eeqar}{\end{eqnarray}}

\newtheorem{proposition}{Proposition}
\newtheorem{leftnumbers}{\rule{-0.35em}{0em}}[section]
\newenvironment{themit}[1]{\begin{leftnumbers}{\sc #1.}}{\end{leftnumbers}}

\newcommand{\btheorem}{\begin{themit}{Theorem}\ }
\newcommand{\etheorem}{\end{themit}}
\newcommand{\bproposition}{\begin{themit}{Proposition}\ }
\newcommand{\eproposition}{\end{themit}}
\newcommand{\blemma}{\begin{themit}{Lemma}\ }
\newcommand{\elemma}{\end{themit}}
\newcommand{\bconjecture}{\begin{themit}{Conjecture}\ }
\newcommand{\econjecture}{\end{themit}}
\newcommand{\bproblem}{\begin{themit}{Problem}\ }
\newcommand{\eproblem}{\end{themit}}
\newcommand{\bcorollary}{\begin{themit}{Corollary}\ }
\newcommand{\ecorollary}{\end{themit}}

\newenvironment{remdefexam}[1]{\begin{leftnumbers}{\sc #1.}
 \rm}{\end{leftnumbers}}

\newcommand{\bremark}{\begin{remdefexam}{Remark}\ }
\newcommand{\eremark}{\end{remdefexam}}
\newcommand{\bdefinition}{\begin{remdefexam}{Definition}}
\newcommand{\edefinition}{\end{remdefexam}}
\newcommand{\bexample}{\begin{remdefexam}{Example}\ }
\newcommand{\eexample}{\end{remdefexam}}
\newcommand{\bexercise}{\begin{remdefexam}{Exercise}\ }
\newcommand{\eexercise}{\end{remdefexam}}

\newcommand{\beqn}{\refstepcounter{leftnumbers}
\begin{equation}}
\newcommand{\eeqn}{\end{equation}}

\newtheorem{remit}[proposition]{Remark}

\newtheorem{defineit}[proposition]{Definition}

\newtheorem{examit}[proposition]{Example}

\newcommand{\cals}{\mbox{$\cal S$}}

\newcommand{\bc}{{\bf C}}

\newcommand{\bz}{{\bf Z}}

\newcommand{\frechet}{Fr\'{e}chet\ }

\newcommand{\bq}{\mbox{\bf Q}}

\newcommand{\tr}{{\rm Tr}}

\title{Chern character for the Schwartz algebra of $p$-adic $GL(n)$}
\author{Jacek Brodzki
\and 
Roger Plymen 
} 
\date{}
\begin{document}
\maketitle

\mysection{Introduction}
Atiyah has demonstrated that rational cohomology of a compact Hausdorff space
can be defined in terms of $K$-theory. 
This is made possible by the existence of a Chern character
$$
ch_j: K^j(X)\ra H^{ev/odd}(X; \bz), \quad j=0,1. 
$$
As in \cite[p.~142]{Atiyah}(see also \cite[Sec.~6.2]{Rosenberg}), this map 
is a rational isomorphism:
\beqn\label{1}
ch_j\otimes_{\bz} 1: K^j(X)\otimes_\bz\bq \simeq 
H^{ev/odd}(X; \bq).
\eeqn

When $X$ is a compact smooth manifold, the value of the Chern character
on a $K$-theory class represented by a bundle $E$ can be computed 
explicitly in terms of an inhomogeneous differential form 
whose components are (up to factors) exterior powers of the curvature
of $E$. Thus, for instance, information about $K$-theory of $X$ may be 
obtained by pairing the cohomology class represented by the 
Chern character with closed de Rham currents. 
There is now a construction, due to Gorokhovsky \cite{Ghor}, 
which provides explicit formulae for the classical Chern character of a compact
Hausdorff space $X$ in terms of Alexander-Spanier cocycles. In the
case when $X$ is a compact smooth manifold, both ways of 
computing the Chern character coincide.

Connes created a Chern character from $K$-homology 
of an algebra $A$, defined in terms of abstract elliptic operators, with 
values in cyclic cohomology \cite[Chap.~I]{Connes}. 
He also provided a key idea for the construction of a dual Chern character 
from the $K$-theory of $A$, defined in terms of analogues of vector 
bundles, with values in cyclic homology \cite[Chap.~II, Prop.~14]{Connes}. 
The two characters
are compatible with the index pairing between $K$-homology and $K$-theory 
of the algebra $A$. 

These groundbreaking ideas were developed by a number of authors 
\cite{JonesKassel,Kassel,Loday,Nistor} and completed by Cuntz and Quillen
who provided explicit formulae for Chern characters with values
in periodic cyclic homology \cite{CQ3}. 
More recently, Cuntz constructed a universal 
bivariant $kk$-functor and a compatible Chern character whose target is the 
bivariant periodic cyclic cohomology of Cuntz and Quillen \cite{Cuntz:Doc}. 

An important application of these constructions is provided by a proof
of the Baum-Connes conjecture for the $p$-adic group $GL(n)$ \cite{BHP:CR}. 
There,  Chern characters were used to obtain information about 
the Baum-Connes assembly map. For a locally compact group $G$
the Baum-Connes map is an index map from the universal example of 
proper $G$ actions to the $K$-theory of the reduced $C^*$-algebra
$C^*_r(G)$. 

In this paper we construct a Chern character
\beqn\label{2}
ch:K_*(C^*_r(G)) \ra HP_*(S(G))
\eeqn
where $S(G)$ is the Schwartz algebra of the $p$-adic group $GL(n)$. 
We then prove that this map is an isomorphism after tensoring 
over $\bz$ with $\bc$. The first step in our construction 
is to adapt the formula of Cuntz and Quillen to the context of \frechet\
algebras. It turns out that this translation is very natural in the case 
when the \frechet\ algebra under study is a dense subalgebra 
of a $C^*$-algebra which is stable under holomorphic functional 
calculus. The Schwartz algebra $S(G)$ is defined as the strict 
inductive limit of nuclear \frechet\ algebras $S(G//K)$ which satisfy 
this property \cite[p.~93]{HC}. 

Using explicit formulae provided by the connection-curvature computation 
of the Chern character in the Chern-Weil theory and Gorokhovsky's results in 
the topological case we compare the map of Cuntz and Quillen with 
the classical Chern character in algebraic topology. Ultimately 
our result that the Chern character (\ref{2}) is an isomorphism after
tensoring with $\bc$ follows from the fact that the classical 
Chern character is a rational isomorphism as stated in (\ref{1}).

\mysection{\frechet algebras and the Chern character}
 
Cuntz and Quillen provide explicit formulae \cite{CQ3} for Chern character maps
$ch_{CQ}: K_*(A) \ra HP_*(A)$ from the $K$-theory of a unital algebra $A$
to its periodic cyclic homology. The goal of this section is to adapt their 
results to the context of \frechet\ algebras and to
show how 
these maps compare, in the case of the algebra of smooth functions on 
a manifold, with the Chern character given by the Chern-Weil theory. 

By a \frechet\ algebra we will understand a locally 
multiplicatively convex
\frechet\ algebra over $\bc$. 
That is, a \frechet\ algebra is a complete metrizable 
topological algebra whose topology is given by a countable 
family of submultiplicative seminorms. 

$K$-theory for \frechet\ algebras 
was constructed by Phillips \cite{Phillips}. 
Any definition  of $K_*(A)$ requires a suitable choice of stabilization 
of $A$. In the case of Banach or $C^*$-algebras one uses 
infinite matrices over $A$ with finitely many nonzero entries, 
that is the direct limit $\bigcup_n M_n(A)$ of the natural 
direct system of matrix algebras $M_n(A)$ of increasing size. 
This approach is not suitable for a general \frechet\ algebra, see 
\cite[Section 2]{Phillips} for a full discussion. 
It turns out that the stabilization of $A$ best adapted
to the definition of $K$-theory for \frechet\ algebras is given 
by $S(\bz^2) \otimes A$, where  $\cals(\bz^2)$ is the algebra
of rapidly decreasing  functions on $\bz^2$. However, there
is a significant simplification in the 
case of algebras which are dense subalgebras of $C^*$-algebras, 
stable under holomorphic functional calculus and which are \frechet\ 
algebras in a finer topology. 
All \frechet\ algebras considered in this paper are 
of this kind unless stated otherwise. In this case Phillips proves that the
$K$-theory  can be defined in the following familiar way
\cite[Def.~7.1, Cor.~7.9]{Phillips}. For a unital \frechet\ algebra $A$,
$K_0(A)$ is defined as the Grothendieck group of the semigroup of isomorphism
classes of finitely generated projective modules. Thus a class in 
$K_0(A)$ can be represented by a projection in $M_n(A)$, 
for some $n$. Moreover, 
projections in the class of a projection  $e$ are homotopic to $e$.  
Let $GL_n(A)$ denote the group of invertible $n\times n$ matrices, and
let $GL_n(A)_0$ be the path  component of the identity. We define
$K_1(A) = \lim GL_n(A)/GL_n(A)_0$. For nonunital algebras $K_i(A)$ 
is defined as the kernel of the map from $K_i(\tilde{A})$ to 
$K_i(\bc)$, where $\tilde{A}$ is the unitization of $A$. 

Cyclic type homology theories of a complex algebra $A$, for example Hochshild, cyclic and 
periodic cyclic homology, are computed using a mixed complex
$(\Omega A, b, B)$ associated with the algebra $A$ \cite{CQ2}. The universal differential 
graded algebra $\Omega A$ is generated by the algebra $A$ and the 
symbols $da$, $a\in A$, which are linear in $A$ and satisfy the Leibniz identity
$d(ab) = d(a) b + adb$ for any $a,b\in A$. In degree $n$, the space of 
$n$-forms $\Omega^n A$ is spanned by the elements $a_0da_1\dots da_n$. 
As a vector space, $\Omega ^n A = A\otimes (A/\bc)^{\otimes n}$. 
The differential 
$b: \Omega^n A\ra \Omega^{n-1}A$ is defined in positive degrees
 by $b(\omega da) = (-1)^{|\omega|}[\omega, a]$, where $|\omega|$ is the degree of 
 the form $\omega$. In degree zero we put $b $ equal identically zero. 
 The degree $+1$ differential $B:\Omega^{n}A \ra \Omega^{n+1} A$ is defined
 by $B = \sum_{i=0}^n \kappa^i d$, where $\kappa$ is the degree zero 
 Karoubi operator: $\kappa (\omega da) = (-1)^{|\omega|} da\cdot \omega$. 
 The two operators $b$ and $B$ anticommute and are of square zero: 
 $$
 b^2 = Bb +bB = B^2 = 0.
 $$
For a locally convex  algebra $A$ one must choose a suitable topological 
tensor product to use in the definition of the differential graded algebra
 $\Omega A$. In this paper we choose the completed inductive tensor product
 $\bar{\otimes}$ as in \cite{BP}, to which we refer the reader for 
 further details. Thus, for all $n\geq 0$, $\Omega^{n}(A) = 
 A\bar{\otimes} (A/\bc)^{\bar{\otimes} n}$. The differentials $b$ and $B$ are
 continuous in this topology. 
 
 Let $A$ be a unital locally convex algebra over $\bc$. By definition, the {\em Hochschild 
 homology} of $A$ is
 $$
 HH_*(A) = H_*(\Omega A, b). 
 $$
The {\em periodic cyclic homology} $HP_*(A)$ of $A$
is  the homology of 
the $\bz/2\bz$-graded complex 
$$
\dots \arr{B-b} \prod_{n\geq 0} \Omega^{2n}(A) \arr{B-b} 
\prod_{n\geq 0} \Omega^{2n+1}(A) \arr{B-b} 
\prod_{n\geq 0}\Omega^{2n}(A) \arr{B-b} \dots
$$

\mbox{}

Let $A$ be a unital \frechet\ algebra. 
Let $ [e]\in K_0(A)$
be the $K$-class of an idempotent matrix $e
\in M_k(A)$ over $A$. Then the even Chern character of Cuntz and 
Quillen assigns to $[e]$ an even class in the periodic cyclic homology 
of $A$ represented by the even periodic cycle
\beqn\label{CQeven}
ch_{CQ}: e \mapsto \tr e + \sum_{n\geq 1} \frac{(2n)!}{n!}
\tr \left(\left((e - \frac{1}{2}\right)d e^{2n}\right) \in \Omega^{ev} A .
\eeqn
Here $\Omega^{ev}A = \prod_{n\geq 0} \Omega ^{2n} A$. 

In the odd case, let us assume that a class $[g]\in K_1(A)$ is represented 
by an invertible matrix $g\in GL_k(A)$. 
The odd Chern character maps the class $[g]$ to the odd periodic 
cyclic homology class represented by the odd cycle
\beqn\label{CQodd}
ch_{CQ}: g\mapsto \sum_{n\geq 0} n! \tr (g^{-1} dg ( dg^{-1} dg)^n)
\eeqn
in $\Omega^{odd}(A) = \prod_{n\geq 0} \Omega^{2n+1} A$. 

This Chern character is compatible with algebra homomorphisms 
in the sense that a homomorphisms $\phi: A\ra B$ of algebras
gives rise to the commutative diagram
$$
\begin{array}{ccc}
K_*(A) &  \arr{ch_{CQ}} & HP_*(A) \\
\dar{\phi} & & \dar{\phi} \\
K_*(B) & \arr{ch_{CQ}} & HP_*(B)
\end{array}
$$
This map  is compatible with Morita equivalence in the following sense. 
There is a similar diagram in which $B$ is replaced
by the matrix algebra $M_n(A)$ and where vertical arrows are isomorphisms. 

\mbox{}

The formula of Cuntz and Quillen works for any \frechet\ algebra, 
commutative or not.  Its main strength lies in the fact that it 
can be adapted to coincide with the classical Chern character
known from algebraic topology or the Chern character defined 
for smooth manifolds using the Chern-Weil theory. 
To make a comparison of the Chern character of Cuntz and Quillen with 
the one known from the Chern-Weil theory, let us recall the main points
in the construction of the latter map.

Let $E$ be a complex vector bundle over a compact smooth manifold 
 $M$ and let $\nabla$ be 
a connection on $E$. The curvature $R$ of $\nabla$ is 
an ${\rm End} E$-valued two-form on $M$. Then $\tr \exp(R/2\pi i)$
is an inhomogeneous even exact form on $M$. The component of degree $2k$ represents
a class in $H^{2k}_{dR}(M)$.  
With this understood, the even Chern-Weil character
is the  map 
\beqn\label{CW}
ch_{CW}: [E]\mapsto \left[\tr \exp(R/2\pi i)\right]\in H^{ev}_{dR}(M)= 
\bigoplus_{n\geq 0}
H^{2n}_{dR}(M), 
\eeqn
where $[E]\in K^0(M)$ is the $K$-theory class represented by the bundle 
$E$. 

In the odd case, we assume that a class in $K^1(M)$ is represented
by an invertible matrix $g\in GL_n(C^\infty(M))$. 
The odd Chern character is defined by the formula \cite[p.491]{Getzler}
\beqn\label{getzler}
ch_{CW}:  [g]\mapsto 
\sum_{n\geq 0} \left[(-1)^n\frac{n!}{(2\pi i)^n(2n+1)!}\tr((g^{-1}dg)^{2n+1})
\right]
\eeqn
which takes values
in  $H^{odd}_{dR}(M)= \bigoplus_{n\geq 0} H^{2n+1}_{dR}(M)$. 
\bremark
The factors $(2\pi i)^{-n}$ in the two formulae above are used in 
algebraic topology to ensure that the Chern characters represent integral 
classes in cohomology of $M$. This plays a crucial role, for instance,
in the index theorems of Atiyah and Singer.
\eremark

When $A$ is the algebra $A=C^\infty(M)$ of smooth 
functions on a smooth compact manifold $M$ the Cuntz-Quillen 
maps (\ref{CQeven}) and (\ref{CQodd}) coincide with 
(\ref{CW}) and (\ref{getzler}) as we now show.  

Let us denote by $(\Omega M, d) $ the complexified de Rham complex 
of the manifold $M$. The universal property \cite{CQ1} 
of the differential graded
algebra $\Omega A$ implies that there is a natural surjection of mixed 
complexes
$$
\mu: (\Omega A , b , B) \ra (\Omega M , 0 , d)
$$
which sends a noncommutative $n$-form $f^0df^1\cdots df^n$ to the 
differential form $(1/n!) f^0df^1\wedge \dots \wedge df^n$. This 
map induces an isomorphism of all cyclic type homology theories
associated with the two mixed complexes if the algebra $A$ satisfies
the Hochschild-Kostant-Rosenberg theorem 
$$
HH_*(A)= H_*(\Omega A, b) = \Omega M.
$$
A crucial result of Connes' shows  that this is true for $A=C^\infty(M)$. 
In particular we have in this case that 
$HP_{ev/odd}(A)\simeq H^{ev/odd}_{dR}(M)$.

The map $\mu$ thus induces a map $\mu: HP_{ev/odd}(A) \ra
H^{ev/odd}_{dR}(M)$. The de Rham cohomology groups are $\bc$-modules
with the natural action $\lambda[\omega] = [\lambda\omega]$ 
for any $\lambda\in \bc$. We use
this fact to define the scaling map
$$
c_{n}: H^n(M) \ra H^n(M)
$$
which for $n=2k$ multiplies by $1/(2\pi i)^k$ and for 
$n=2k+1$ by $1/(2 \pi i)^{k+1}$. 
We let $c: \bigoplus_n H^n(M)\ra \bigoplus H^n(M)$ 
be the map whose component in degree $n$ is $c_n$. It is clear 
that $c$ is an isomorphism of $\bc$-modules. Let $\mu_c= c\circ \mu$. 

The comparison of the two Chern characters is established by assembling
various known results which we recall here for future reference. 
\bproposition\label{CQsmooth}
Let $M$ be a smooth compact manifold and let  $A= C^\infty(M)$
be the algebra of smooth functions on $M$. 
The following diagram commutes
$$
\begin{array}{ccc}
K_*(A) & \arr{ch_{CQ}} & HP_*(A) \\
\dar{SS} & & \dar{\mu_c}\\
K^*(M) & \arr{ch_{CW}} & H^{ev/odd}_{dR}(M)
\end{array}
$$
where the  left vertical arrow is the isomorphism given by the
Serre-Swan theorem.  
\eproposition
\proof 
In the even case, the Serre-Swan map sends a projection 
$e\in M_k(C^\infty(M))$ 
to the bundle $E=e\theta^k$, where $\theta^k$ is the trivial 
complex rank-$k$ bundle on $M$. The trivial bundle $\theta^k$ 
is equipped with the canonical connection $d$, which is 
just the de Rham differential.  
The subbundle $E$  of $\theta^k$ 
is then equipped with the induced Grassmannnian 
connection 
$$
D = e \cdot d \cdot e
$$
which acts on sections of $E$. A typical section of $E$ is of the 
form $ef$,  where $f\in C^\infty(M, \bc^k)$ is a section of the trivial 
bundle $\theta^k$. The curvature of this connection is (cf. \cite[p.223]{Quillen:Cayley})
$
R = e(de)^2.
$
Moreover,  $e(de)^2 = (de)^2e$ so that 
$
R^k = e(de)^{2k}.
$
Hence the Chern character of the sub-bundle $E$ is 
$$
ch(E) = \exp(R/2\pi i) = \sum_{n\geq 0} \frac{1}{(2\pi i)^nn!} e(de)^{2n}
$$
where the power on the right is taken using the exterior product of 
forms. 

The commutativity of the diagram in the even case is now established 
using that \cite[p.~435]{CQ3}
$$
\mu_c ch_{CQ}([e]) = \sum_{n\geq 0}\frac{1}{(2\pi i)^n n!}\tr (e(de)^{2n})
= ch_{CW}(SS[e])
$$
The odd case follows once we notice that
\cite[p.~436]{CQ3}
$$
\mu_c[g]= \sum_{n\geq 0} \frac{n!}{(2\pi i)^n(2n+1)!}\tr ((g^{-1}dg)^{2n+1}).
$$
\qed

\mysection{The Schwartz algebra of  $GL(n)$}

Let $F$ be a nonarchimedean local field and let $G= GL(n)=GL(n,F)$. 
Let $K$ be a compact open subgroup of $G$. We define $S(G//K)$ 
to be the algebra of all complex-valued
 functions on $G$ which 
are rapidly decreasing and $K$-bi-invariant,
 with product given by convolution.
 \blemma\label{Mischenko}
For every compact open subgroup $K$ of $G$ the 
algebra $S(G//K)$ is an $m$-convex  \frechet\ algebra. 
\elemma
\proof 
Mischenko's Fourier transform \cite{Mischenko} provides the following
isomorphism of locally convex algebras
\beqn\label{MiscFormula}
S(G//K)  \simeq \displaystyle{\bigoplus_M} [C^\infty({\rm End}
 F(M:K))]^{W(M)} 
\eeqn
where $F(M:K)\rightarrow X(M:K) $ is a complex, Hermitian, trivialized 
$W(M)$-bundle over a compact manifold $X(M:K)$. The group $W(M)$ is 
a finite reflection group and $M$ is a Levi subgroup. One Levi subgroup
is chosen in each conjugacy class of $G$; in particular, 
the sum on the right is finite. 
Each direct summand on the right is a subalgebra of the algebra
$C^\infty(X(M:K)) \otimes M_n(\bc)$ of  matrix-valued smooth
functions on $X(M:K)$. 
The algebra $C^\infty(V)$ of smooth functions on a compact smooth 
manifold $V$ is an $m$-convex \frechet\ algebra.
Indeed, for every natural number $n$ we 
put $p_n(f) = \sup \{|\alpha|\leq n, \, x\in V \mid 
|\partial^\alpha f(x)|\}$ and define
$$
q_n(f) = \sum_{i=0}^n \frac{1}{i!} p_i(f)
$$
The seminorms $q_n$ are submultiplicative. 

The matrix 
algebra $M_n(\bc)$ is  equipped with the Hilbert-Schmidt norm $\norm{\; }_2$, 
which is submultiplicative. 

The algebraic tensor product $C^\infty(X(M:K))\otimes M_n(\bc)$ 
is topologised by 
submultiplicative seminorms $q_n\otimes \norm{\;}_2$ and so is an $m$-convex
\frechet\ algebra. 
This remark combined with Mischenko's isomorphism shows that 
$S(G//K)$ is an $m$-convex \frechet\ algebra for every compact open 
subgroup $K$ of $G$. 
\qed

Following Harish-Chandra \cite[p.~93]{HC}, the 
Schwartz algebra of $G$ is defined as the strict inductive limit 
of the  algebras $S(G//K)$, that is
$$
S(G) = \lim_{\arr{}} S(G//K) = \bigcup_K S(G//K).
$$
Choose a left-invariant Haar measure on $G$. Then $L^1(G)$ acts on $L^2(G)$ by 
convolution
$$
\lambda(f)h = f*h, 
$$
for $f\in L^1(G)$ and $h\in L^2(G)$. The reduced 
$C^*$-algebra $C^*_r(G)$ is the closure of $\lambda(L^1(G))$ in 
the $C^*$-algebra of bounded operators on $L^2(G)$. 

In this section, we  prove the following main result of this paper. 
\btheorem\label{OurResult}
There exists a Chern character map
\beqn\label{Main}
ch: K_*(C^*_r(G)) \ra HP_*(S(G))
\eeqn
which is an isomorphism after tensoring over $\bz$ with $\bc$. 
\etheorem

\proof
For every compact open subgroup $K$ of $G$ the algebra $S(G//K)$ 
is stable under holomorphic functional calculus, as follows from 
Mischenko's theorem (\ref{MiscFormula}). Moreover,  each algebra
$S(G//K)$ is an $m$-convex \frechet\ algebra by Lemma
\ref{Mischenko}.  
Denote by $C^*_r(G//K)$ the closure of $S(G//K)$  in the 
reduced $C^*$ algebra $C^*_r(G)$. Then $S(G//K)$ is 
dense in $C^*_r(G//K)$. Thus the Cuntz-Quillen formula applies to 
give a Chern character
$$
ch^K_{CQ}: K_*(S(G//K)) \ra HP_*(S(G//K))
$$
This  Chern character then determines a  map: 
$$
ch^K_{CQ}: K_*(C^*_r(G//K)) \simeq K_*(S(G//K)) \ra HP_*(G//K)
$$
The $C^*$-algebras $C^*_r(G//K)$ form a direct system whose
$C^*$-inductive limit is the reduced $C^*$-algebra $C^*_r(G)$. 
We thus have a direct system of Chern characters
$ch^K_{CQ}$ and we define 
$$
ch = \displaystyle{\lim_{\arr{}}}\; ch^K_{CQ}: K_*(C^*_r(G)) \ra 
HP_*(S(G))
$$ 
We use here continuity of $K$-theory and periodic cyclic homology with 
respect to direct limits which holds in this case \cite[Theorem 6]{BP}. 

We now prove that this character is an  isomorphism after tensoring 
(over $\bz$) with $\bc$.  

Relying on  Mischenko's theorem stated in Lemma \ref{Mischenko} we proved
in \cite{BP} that each algebra $S(G//K)$ is Morita equivalent to a commutative
algebra
\beqn\label{Morita}
S(G//K) \stackrel{\rm Morita}{\sim} 
\displaystyle{\bigoplus_M} C^\infty(X(M:K))^{W(M)}
\eeqn
where each direct summand is the algebra of invariant smooth functions on 
$X(M:K)$.

Using this fact together with the compatibility of functors $K_*$ and $HP_*$ 
with algebra homomorphisms and Morita equivalence we see that, for 
every open subgroup  $K$ of $G$, 
there is the commutative diagram
\beqn\label{Levi}
\begin{array}{ccc}
K_*(S(G//K)) & \arr{ch^K_{CQ}} & HP_*(S(G//K)) \\
\Bdal{\simeq} && \Bdal{\simeq}\\
\bigoplus_M K_*(C^\infty(X(M:K)))^{W(M)} & \arr{ch_{CQ}}& \bigoplus_M
HP_*(C^\infty(X(M:K)^{W(M)}))
\end{array}
\eeqn
Vertical arrows in this diagram  are isomorphisms. The bottom arrow is 
the direct sum of Cuntz-Quillen Chern characters defined for 
each Levi subgroup $M$ using formulae (\ref{CQeven}) and 
(\ref{CQodd}). The key 
step in our argument is to show that each of the components of the
map at the bottom is an isomorphism after tensoring with $\bc$. 

This follows from the following  geometric consideration. 
Let $X$ be a compact $W$-manifold, where $W$ is a finite reflection group. 
There is the following commutative diagram
\beqn\label{triplediagram}
\begin{array}{ccc}
K_*(C^\infty(X)) & \arr{ch_{CQ}} & HP_*(C^\infty(X))\\
\parallel && \dar{\mu_c} \\
K_*(C^\infty(X)) & \arr{ch_{CW}} & H^*_{dR}(X)\\
\parallel &&\uar{\lambda}\\
K_*(C^\infty(X)) & \arr{ch_G} & H^*_{AS}(X).
\end{array}
\eeqn
The top two rows form the commutative diagram of Proposition 
\ref{CQsmooth}. The bottom row introduces the Chern character
with values in the Alexander-Spanier cohomology $H^*_{AS}(X)$, 
which was
constructed  by Gorokhovsky for any 
compact Hausdorff space
\cite{Ghor}. Ghorokhovsky proves that his map coincides with the 
classical Chern character. In particular, it is a 
rational isomorphism by \cite{Karoubi:CartanSem}. 
Moreover, when $X$ is a compact smooth manifold, there is a  
natural isomorphism 
$$\lambda : H^*_{AS}(X) 
\ra H^*_{dR}(X)$$
which identifies the Alexander-Spanier cohomology with the 
de Rham cohomogy. It is shown in \cite{Ghor}
that this  map sends Gorokhovsky's Chern character 
to the Chern character known from the Chern-Weil theory. 
In other words, the bottom part of the above diagram commutes. 

It follows from the functorial properties of the theories 
involved that this diagram is compatible with diffeomorphisms 
of $X$; in particular all maps  are 
$W$-maps. 

Let us now consider an idempotent $e\in M_k(C^\infty(X)^W)$. 
For any $w\in W$ we have that $w^*[e] = [w^*e]=[e]$. For, if 
$f$ is path-homotopic to $e$ via a path $e_t$ of idempotents, 
then $w^*f$ is connected to $w^*e=e$ by the continuous path $w^*e_t$. 
Thus  the $K$-class of $e$ is invariant under
the action of $W$. In the odd case, 
the class of an invertible matrix $g\in GL_k(C^\infty(X)^W)$ constist
of elements of this group which are homotopic to $g$. But the 
action of $W$ sends invertible elements to invertible elements
and if $h$ is connected to $g$ by a continuous path of invertibles, 
then the same is true about $w^*g = g$ and $w^*h$. 

By functoriality, the image of $K_*(C^\infty(X)^W)$ under the Chern character 
map in each of the three cases  must therefore 
be a submodule of the invariant part of the corresponding cohomology theory on 
the right. In other words, we have the following commutative diagram. 
$$
\begin{array}{ccc}
K_*(C^\infty(X)^W) & \arr{ch_{CQ}} & HP_*(C^\infty(X))^W\\
\parallel && \dar{\mu_c} \\
K_*(C^\infty(X)^W) & \arr{ch_{CW}} & H^*_{dR}(X)^W\\
\parallel &&\uar{\lambda}\\
K_*(C^\infty(X)^W) & \arr{ch_G} & H^*_{AS}(X)^W
\end{array}
$$
The two vertical arrows on the right remain isomorphisms, since the
original isomorphisms were compatible with the action of the 
group $W$. 

Since $C^\infty (X)^W$ holomorphically closed and dense in $C(X)^W$
the $K$-theory of the two algebras is the same
$$
K_*(C^\infty(X)^W) \simeq K_*(C(X)^W). 
$$
Furthermore, the isomorphism of $C^*$-algebras $C(X)^W= C(X/W)$ gives
$$K_*(C(X)^W) = K_*(C(X/W) =K^*(X/W),$$ 
where the last identity follows from
the Serre-Swan theorem using the fact that $X/W$ is a compact Hausdorff 
space. Thus we have established the isomorphism
$$
K_*(C^\infty(X/W)) =  K^*(X/W).
$$
On the side of homology it is well known that  $H_{AS}^*(X)^W = H^*_{AS}(X/W)$. 
On the other hand the result of Wassermann 
\cite[p.~238]{Wassermann} gives that 
$HP_*(C^\infty(X)^W) = HP_*(C^\infty(X))^W$. 
We thus have the following
commutative diagram
$$
\begin{array}{ccc}
K_*(C^\infty(X)^W) & \arr{ch_{CQ}} & HP_*(C^\infty(X)^W)\\
\Bdal{\simeq} && \Bdal{\lambda^{-1}\circ \mu_c} \\
K^*(X/W) & \arr{ch_G}& H^*_{AS}(X/W,\bc)
\end{array}
$$
in which the bottom map is the Chern character constructed 
by Gorokhovsky. Since Gorokhovsky's Chern character is 
the same as the classical Chern character, the bottom 
map becomes an isomorphism after the left hand side
has been tensored with $\bc$.
Applying this 
reasoning to the situation described by the diagram 
\ref{Levi} we see that the bottom arrow in that diagram 
is an isomorphism after tensoring with $\bc$.

We have therefore proved that for every compact open subgroup $K$ the map 
$$
ch^K_{CQ}: K_*(C^*_r(G//K))\otimes_{\bz}\bc \ra HP_*(S(G//K))
$$
is an isomorphism. Passing to the limit we finish the proof of the Theorem. 
\qed

\bremark
The result stated in Theorem \ref{OurResult} is crucial
 to the proof of the Baum-Connes conjecture given 
in \cite{BHP:CR}. 
\eremark

\noindent
{\footnotesize J.B.: School of Mathematical Sciences, University of Exeter, Exeter EX4 4QE;
brodzki@maths.ex.ac.uk\\
R.P.: Department of Mathematics, University of Manchester, Manchester M13 9PL;
roger@ma.man.ac.uk}


\begin{thebibliography}{9999}
\bibitem{Atiyah} M.~F.~Atiyah, $K$-theory, Benjamin, New York, 1967. 
\bibitem{BHP:CR} P.~Baum, N.~Higson, R.~J.~Plymen, A proof of the 
Baum-Connes conjecture for $p$-adic $GL(n)$,  
{\em C.~R. Acad. Sci. Paris} {\bf 325} (1997) 171-176. 
\bibitem{BP} J.~Brodzki, R.~J.~Plymen,  Periodic cyclic homology 
of certain nuclear algebras,  {\em C. R. Acad. Sci. Paris}
 {\bf 329} (1999), 671-676.
\bibitem{Connes} A.~Connes, Noncommutative differential geometry, 
{\em Publ.~Mathem. I.H.E.S.} {\bf 62} (1985), 257-360. 
\bibitem{Connes:Book} A.~Connes, Noncommutative geometry, Academic 
Press, London, New York 1994. 
\bibitem{CM} A.~Connes, H.~Moscovici,  Cyclic cohomology, the Novikov
conjecture and hyperbolic groups,  {\em Topology} {\bf 29}, (1990), 
345-388. 
\bibitem{Cuntz:Doc} J.~Cuntz, Bivariante $K$-theorie f\"{u}r lokcalconvexe
Algebren und der Chern-Connes-Charakter, {\em Doc. Math} {\bf 2} (1997), 
139-182.
\bibitem{CQ1} J.~Cuntz, D.~Quillen,  Algebra extensions and nonsingularity, {\em Journal AMS} {\bf 8} (1995), 251-289. 
\bibitem{CQ2} J.~Cuntz, D.~Quillen, Operators on noncommutative differential forms 
and cyclic homology, In: Geometry, Topology and Physics for Raoul Bott, International Press, 
Cambridge MA, 1995. 
\bibitem{CQ3} J.~Cuntz, D.~Quillen, Cyclic homology and nonsingularity, 
{\em Journal AMS}, {\bf 8}, 373-442. 
\bibitem{Getzler} E. Getzler, 
The odd Chern character in cyclic homology and spectral flow, 
{\em Topology} {\bf 33} (1994), 663-681. 
\bibitem{Ghor} A. Gorokhovsky, Chern classes in Alexander-Spanier Cohomology, {\em $K$-theory}
{\bf 15} (1998), 253-268.
\bibitem{HC} Harish-Chandra, Collected Papers, Vol.~4, Springer Verlag, 
Berlin 1984. 
%\bibitem{Hirzebruch} F.~Hirzebruch, Topological methods in algebraic
%geometry, Springer Verlag, Heidelberg 1978. 
\bibitem{JonesKassel} J.D.S.~Jones, C.~Kassel, Bivariant cyclic 
theory, {\em $K$-theory} {\bf 3} (1989), 339-366. 
\bibitem{Kassel} C.~Kassel, Caract\'{e}re de Chern bivariant, 
{\em $K$-theory} {\bf 3} (1989), 367-400. 
%\bibitem{Karoubi:Ktheory} M.~Karoubi, $K$-theory, an introduction, 
%Springer Verlag, Heidelberg, New York 1978. 
\bibitem{Karoubi:CartanSem} M.~Karoubi, Le isomorphismes
de Chern et de Thom-Ghysin en $K$-th\'{e}orie, 
 S\'{e}minaire  Cartan-Schwartz, 
16e ann\'{e}e, 1963/64, n. 16. 
\bibitem{Loday} J.-L.~Loday , Cyclic homology, Springer-Verlag, 
Berlin, 1992.
\bibitem{Mischenko} P.~Mischenko, 
Invariant tempered distributions on the reductive 
$p$-adic group $GL_n(F_p)$, C.~R.~Math. Rep. Acad. Sci. Canada 4 (1982) 123-127
\bibitem{Nistor} V.~Nistor, {\em A bivariant Chern-Connes character}, 
Annals of Math. {\bf 138} (1993), 555-590. 
\bibitem{Phillips} N.C. Phillips, $K$-theory for \frechet\ algebras, 
{\em Int. J. Mathem.} {\bf 2} (1991), 77-129. 
\bibitem{Quillen:Cayley} D.~Quillen, Superconnection
character forms and the Cayley transform, {\em Topology} {\bf 27} (1988), 
211-238. 
\bibitem{Rosenberg} J.~Rosenberg, Algebraic $K$-theory and its applications, 
GTM 147, Springer Verlag, Berlin, New York 1994. 
\bibitem{Wassermann} A.~Wassermann,  Cyclic cohomology of algebras 
of smooth functions on
orbifolds, London Math.~Soc.~Lecture Notes 135 (1988) 229-244
\end{thebibliography}
\end{document}